%% file: main.tex
\documentclass{article}
\usepackage{spconf,amsmath,graphicx}
\usepackage{amsmath}
\usepackage{mathtools}
\usepackage{amsfonts}
\usepackage{commath}
\usepackage{graphicx}    
\usepackage{enumitem}
\usepackage{tcolorbox}
\usepackage{caption}
\usepackage{pifont}
\usepackage{subcaption}

\def\L{{\cal L}}
\newcommand{\Sc}{\ensuremath{\mathcal{S}}}

\newcommand{\eR}{\ensuremath{\mathbb{R}}}
\newcommand{\eN}{\ensuremath{\mathbb{N}}}
\newcommand{\eS}{\ensuremath{\mathbb{S}}}
\newcommand{\xmark}{\ding{55}}

\newtheorem{theorem}{Theorem}
\newtheorem{assumption}{Assumption}

\title{BLOCK DISTRIBUTED 3MG ALGORITHM AND ITS APPLICATION TO 3D IMAGE RESTORATION}
\name{\text{Mathieu Chalvidal}, \text{Emilie Chouzenoux}\thanks{This work was partly supported by the French Agence Nationale de la Recherche under grant ANR-17-CE40-0004-01 MAJIC and by the European Research
Council Starting Grant MAJORIS ERC-2019-STG-850925.}}
\address{Université Paris-Saclay, CentraleSupélec, Inria, CVN, Gif-sur-Yvette, 91190, France}
\begin{document}
%
\maketitle
\ninept
\begin{abstract}
Modern 3D image recovery problems require powerful optimization frameworks to handle high dimensionality while providing reliable numerical solutions in a reasonable time. In this perspective, asynchronous parallel optimization algorithms have received an increasing attention by overcoming memory limitation issues and communication bottlenecks. In this work, we propose a block distributed Majorize-Minorize Memory Gradient (BD3MG) optimization algorithm for solving large scale non-convex differentiable  optimization problems. Assuming a distributed memory environment, the algorithm casts the efficient 3MG scheme~\cite{Chouzenoux_2013} into smaller dimension subproblems where blocks of variables are addressed in an asynchronous manner. Convergence of the sequence built by the proposed BD3MG method is established under mild assumptions. Application to the restoration of 3D images degraded by a depth-variant blur shows that our method yields significant computational time reduction compared to several synchronous and asynchronous competitors, while exhibiting great scalability potential.
\end{abstract}
\begin{keywords}
Majorization-Minimization ; Block-alternating optimization ; Distributed scheme ; Asynchronous communication ; Image deblurring ; Depth-varying blur.
\end{keywords}
\section{Introduction}
\label{sec:intro}

Constantly improving image acquisition devices, from microscopes to medical imaging machines, impose to work with increasingly large data. From the mathematical perspective, many problems of image processing require to solve
\begin{equation}\label{eq:P}
    \operatorname{minimize} \limits_{x\in \mathbb{R}^{N}}f(x)    
\end{equation}
with $f: \mathbb{R}^N \mapsto \mathbb{R}$ a differentiable objective function. To limit the dependence of the optimization process on the dimension $N$ of the problem, block alternating algorithms have been developed. In these schemes, at each iteration only a subset of the variables are updated, by minimizing $f$ with respect to only those variables, the others being fixed. The blocks are selected iteratively following a cyclic (or quasi-cyclic) order or a random rule. For a majority of problems encountered in image restoration, the minimization of $f$ with respect to a given block of variables is usually not possible in a closed form. Furthermore, the application of such basic block coordinate descent update is usually non desirable, as it may lead to convergence issues~\cite{Tseng2001}. Better performance and stability are obtained when combining the block alternating approach with a so-called majoration-minimization (MM) scheme \cite{4303144}. It consists in building, at each iteration, a majorizing approximation for $f$ within the active block of variables, whose minimizer has a more tractable form. Many powerful algorithms fall within this framework, such as BSUM \cite{Hong2016}, PALM \cite{Bolte2014}, NMF \cite{LeeNMF}, to name a few. By relying more closely on the structure of the objective function, block alternating MM methods can reach fast convergence rate \cite{Fessler97,Repetti2015} while offering theoretical guarantees in non-convex cases \cite{Bolte2014,Chouzenoux15jogo}.

When the problem size becomes increasingly large, running such algorithm becomes difficult, due to memory limitation issues. Parallel implementations of MM schemes have been devised, where the block updates are performed simultaneously, allowing to distribute computations on different nodes (or machines)~\cite{Fessler2002,7532949,Moon2019}. Implementation on parallel architecture requires to pay attention to communication cost. The latter can be reduced by resorting to an asynchronous parallel implementation. Each computation node has its own iteration loop, so that it can keep updating its local variables without having to wait for the update of the other, distant variables. This raises challenging questions, in terms of convergence analysis, as the communication delays may introduce instabilities. A plethora of recent works have focused on proposing distributed optimization algorithms with assessed convergence, based on stochastic proximal primal \cite{niu2011hogwild,grishchenko2018asynchronous} or primal-dual \cite{Zhang:2014:ADA:3044805.3045082,pesquet2014class,Hannah_2017,Abboud2019,Onose16,Xu2019} techniques. However, as they rely on the formulation of dual instances of a stochastic coordinate descent strategy, those algorithms are limited to convex (sometimes even strongly convex) optimization and often require specific probabilistic assumptions on the block update rule difficult to meet in practice. In the context of MM algorithms, although the need for distributed implementation strategies is crucial (see the discussion in \cite{Hong2016}), few results are available so far regarding theoretical convergence guarantees. For example, we can mention the work of~\cite{DavisNIPS,davis2016asynchronous}, that proposes an asynchronous version of PALM, with proven convergence in non-convex case, and good practical behaviour~\cite{Thouvenin}. 


%
%
 In this paper, we focus on the MM Memory Gradient (3MG) method \cite{Chouzenoux_2013}. This algorithm integrates a subspace acceleration strategy in the MM framework, leading to one of the most efficient strategies for smooth optimization at large scales \cite{Sun2016}. Here, we propose to set forth an asynchronous distributed version of this algorithm, called Block Distributed 3MG (BD3MG). The latter allows to account efficiently for implementation on parallel architecture with communication delays. By relying on the theoretical framework introduced in \cite{davis2016asynchronous}, we derive convergence guarantees for BD3MG in the challenging non-convex setting, under mild assumptions on the block update rule. We illustrate its performance on the problem of 3D image restoration in the presence of depth-variant blur. The rest of the paper is organized as follows. In Section 2, we introduce our algorithm, analyse its principle features and state its convergence properties. Section 3 focuses on its application to 3D image deblurring. We conclude in Section 4.

\section{Proposed method}
\subsection{Notations}
\label{sec:pagestyle}

We consider the Hilbert space $\mathbb{R}^N$ endowed with the usual scalar product and norm $(\langle \cdot,\cdot \rangle,\| \cdot \|)$. Let a subset \smash{$\Sc \subset [\![1,N]\!]$} of cardinality \smash{$|\Sc|$}. We denote $x_{(\Sc)} = (x_i)_{i \in \Sc} \in \mathbb{R}^{|\mathcal{S}|}$ the restriction of $x$ to the coordinates in $\Sc$. In the same fashion, $\nabla_{(\mathcal{S})}f(x) \in \mathbb{R}^{|\mathcal{S}|}$ denotes the gradient of $f$, evaluated at $x$ along the components indexed in $\mathcal{S}$. Finally, for a matrix $\smash{\mathcal{A} \in \mathbb{R}^{N \times N}}$, we define the submatrix $\mathcal{A}_{(\mathcal{S})} = \left([\mathcal{A}]_{p,p} \right)_{p \in \mathcal{S}} \in\mathbb{R}^{|\mathcal{S}| \times|\mathcal{S}|}$.

\subsection{Block quadratic majorant function}
\label{sec:Blockquadmaj}

Our approach relies on the use of block quadratic majorant functions which constitute quadratic surrogates functions for the restriction of $f$ to any set of coordinates $\mathcal{S} \subset [\![1,N]\!]$. Let $\widetilde{x} \in \mathbb{R}^N$. Let us define :
\begin{multline}
\forall \; v \in \mathbb{R}^{|\mathcal{S}|}, \;\; Q_{(\mathcal{S})}(v,\widetilde{x}) = f(\widetilde{x}) + \langle \nabla_{(S)}f(\widetilde{x}), v-\widetilde{x}_{(S)} \rangle \\
 + \frac{1}{2} \langle v-\widetilde{x}_{(S)}, \mathcal{A}_{(S)}(\widetilde{x})(v-\widetilde{x}_{(S)}) \rangle.
\label{eq:maj}
\end{multline}
Matrix $\mathcal{A}_{(S)}(\widetilde{x}) \in \mathbb{R}^{|\mathcal{S}| \times |\mathcal{S}|}$ hereabove is a symmetric definite positive matrix whose expression depends on $\widetilde{x}$. Following the MM paradigm \cite{4303144}, it should be chosen so as to fulfill: 
\begin{equation}\label{Eq_PropMaj}
\forall v \in \mathbb{R}^{|\mathcal{S}|}, \quad Q_{(\mathcal{S})}(v,\widetilde{x}) \geq f_{(\mathcal{S})}(v;\widetilde{x}),
\end{equation}
where $v \mapsto f_{(\mathcal{S})}(v;\widetilde{x})$ denotes the restriction of function $f$ to coordinates in $\mathcal{S}$, the other coordinates being fixed to those in vector $\widetilde{x}$. The existence and construction for such majorant matrices, not discussed here due to the lack of space, is addressed for instance in \cite{7532949,Chouzenoux15jogo,Sun2016}.

\subsection{Block distributed MM memory gradient algorithm}
\label{sec:majhead}
We are now ready to present our BD3MG algorithm, assuming a memory-distributed \textit{star} cluster of C computing agents with a Master node connected to all other agents. The Master loop starts at $x^0 \in \mathbb{R}^N$ and generates the sequence of iterates $(x^k)_{k\in \mathbb{N}}$, that is incremented whenever a worker $c \in [\![1,C]\!]$ updates a subset of coordinates of the global variable. The main particularity of our algorithm is that the updates can occur in an asynchronous fashion, thus reducing considerably idle time. We denote $\Sc_c^k \subset [\![1,N]\!]$ the processing set of each worker $c$ at times $k$, and $\eS_k = \bigcup_{1 \leq c \leq C} (\Sc_c^k)$ the total set of active blocks at that time. The subset associated to any worker is allowed to change from one iteration to an other. We only impose that there is no overlap in the coordinates updated by the workers at a given time:
\begin{equation}
\label{eq:disj}
(\forall k \in \eN)\quad \bigcap\limits_{c \in \{1,\ldots,C\}}\mathcal{S}_c^k = \emptyset.
\end{equation}
At iteration $k\in \eN$, the Master receives an updated increment \smash{$d_{(\mathcal{S}_c^k)}$} from a given worker $c \in [\![1,C]\!]$. The later is used to increment the corresponding indexes $\mathcal{S}_c^k$ within the global variable $x^{k-1}$, while the others remain untouched, which defines $x^{k}$. The Master then decides for the new set of indexes $\Sc_c^{k+1}$ to be treated by worker $c$. He informs worker $c$ of his new task, and send him the triplet\footnote{Note that communication cost can be reduced, when function $f$ reads as composition of terms with sparse operators. See our discussion in \cite{7532949} for an example.} $(x^{k}, \mathcal{S}_c^{k+1}, (x^{k}-x^{k-1}))_{(\Sc_c^{k+1})}$. 

From the viewpoint of the workers, for each new triplet $(\Sc,x,d_{(\Sc)})$ sent by the Master, a 3MG iterate is performed. The later corresponds to the minimization of the block quadratic majorant function \eqref{eq:maj}, within the memory gradient subspace spanned by the two columns of $
    \mathcal{D}_{(\mathcal{S})}(x) = 
    [-\nabla_{(\mathcal{S})} f(x)\; | \; d_{(\mathcal{S})}] \in \mathbb{R}^{|\mathcal{S}|\times2}.
$
This amounts to find $\widehat{u}$, a minimizer of $u \to Q_{(\mathcal{S})}(x + \mathcal{D}_{(\mathcal{S})}(x)u,x)$, which can be obtained through:
\begin{equation}
\label{eq:minimizer}
   \widehat{u} =  (\mathcal{D}_{(\mathcal{S})}^{\top}(x)\mathcal{A}_{(\mathcal{S})}(x)\mathcal{D}_{(\mathcal{S})}(x))^{\dagger} \mathcal{D}_{(\mathcal{S})}(x)\nabla_{(\mathcal{S})}f(x),
\end{equation}
with $\cdot^\dagger$ the Moore-Penrose pseudo-inverse operator.
The upcoming tables summarize our BD3MG algorithm, that consists of two parts, one to be executed by a 'Master' computing node that receives updates and sends tasks, a second one to be executed by all other computing nodes.
\begin{tcolorbox}[title=Block Distributed 3MG (Master)]
    \label{alg:BMMD_master}
    \small
    $\begin{cases}
            \textbf{Initialization} :\\
            \text{Set } k=0, \;\; x^0 \in \mathbb{R}^N. \;\; \\
            \text{For all } c \in[\![1,C]\!],
            \text{      set } \mathcal{S}_c^0\subset [\![1,N]\!] \text{ s.t.} \bigcap\limits_{c \in[\![1,C]\!]}\mathcal{S}_c^0= \emptyset, \\
            \text{ and send } (x^0,\mathcal{S}_c^0,0_{|\mathcal{S}_c^0|}) \text{ to worker } c.\\
            \text{Define }\eS_0 = \bigcup_{c \in[\![1,C]\!]} \Sc_c^0. \\
            \textbf{While}\text{ a stopping criterion is not met:}\\
            \text{      }\small{(0)} \text{ Wait for any worker to send an update }  \\
            \text{      }\small{(1)} \text{ Receive } (d_{(\mathcal{S}_c^k)}) \text{ from a worker } c \\
            \\
            \text{      }\small{(2)} \text{ Update } \begin{cases}
            x^{k+1}_{(\mathcal{S}_c^k)} = x^{k}_{(\mathcal{S}_c^k)} + d_{(\mathcal{S}_c^k)}\\
            x^{k+1}_{(\overline{\mathcal{S}_c^k})} = x^{k}_{(\overline{\mathcal{S}_c^k})}
            \end{cases}
            \\ 
            \\
            \text{      }\small{(3)} \text{ Choose } \mathcal{S}_{c}^{k+1} \subset[\![1,N]\!] / (\eS_k / \Sc_c^k) \\
                        \qquad \text{For every } c' \in [\![1,C]\!]/\{c\} \text{, set } \Sc_{c'}^{k+1} = \Sc_{c'}^{k}. \\
             \qquad \text{Define } \eS_{k+1} = (\eS_k/ \mathcal{S}_{c}^{k}) \cup \mathcal{S}_{c}^{k+1}\\
            \\
            \text{      }\small{(4)} \text{ Send } (x^{k+1},\mathcal{S}_{c}^{k+1} ,(x^{k+1}-x^{k})_{(\mathcal{S}_{c}^{k+1} )}) \text{ to worker } c \\
            \text{      }\small{(5)} \text{ Increment } k =  k+1 \\
           \end{cases}$
\end{tcolorbox}
\begin{tcolorbox}[title=Block Distributed 3MG (Worker)]
    \label{alg:BMMD_worker}
  $\begin{cases}
        \textbf{While}\text{ the Master stopping criterion is not met:}\\
        \text{      }\small{(1)} \text{ Receive } (x, \mathcal{S}, d_{(\mathcal{S})}) \text{ from Master} \\
        \text{      }\small{(2)} \text{ Set } \mathcal{D}_{(\mathcal{S})}(x) = [-\nabla_{(S)} f(x)\;|\;d_{(\mathcal{S})}] \\
        \text{      }\small{(3)} \text{ Compute } \mathcal{A}_{(\mathcal{S})}(x) \text{ and } \nabla_{(S)} f(x)\\
        \text{      }\small{(4)} \text{ } \mathcal{B}_{(\mathcal{S})}(x) =  (\mathcal{D}_{(\mathcal{S})}(x)\mathcal{A}_{(S)}(x)\mathcal{D}_{(\mathcal{S})}(x))^{\dagger}\\ 
        \text{      }\small{(6)} \text{ } d'_{(\mathcal{S})} = - \mathcal{D}_{(\mathcal{S})}(x)\mathcal{B}_{(\mathcal{S})}(x)\mathcal{D}_{(\mathcal{S})}(x)\nabla_{(S)} f(x)\\
        \text{      }\small{(7)} \text{ Send } (d'_{(\mathcal{S})}) \text{ to the Master} 
    \end{cases}$
\normalsize
\end{tcolorbox}
\subsection{Analysis}
\label{ssec:analysis}
In contrast with its parallel variant \cite{7532949}, BD3MG does not impose any locking condition between workers to perform their computations. Therefore, latency may appear in local variables. A local coordinate $x_i$ used by any worker to perform its update at time $k$ belongs to a previous element \smash{$x^{k'_i}_i$} of the sequence $\{x^k\}_{k \in \mathbb{N}}$ with: $k'_i = \max\{k'\in [\![0,k]\!] \;|\; i \in \mathbb{S}_{k'}\}$. We will model this latency at any iteration $k \in \mathbb{N}$ by introducing $\delta_{k,n}=k-k'_n\in [\![0,k]\!]$ the delay at a coordinate $n \in [\![1,N]\!]$ and $\delta_k = (\delta_{k,n})_{n \in  [\![1,N]\!]}$ the complete vector of delays. For the sake of readability, we will denote $x^{k,\delta_k} = \smash{(x_n^{k-\delta_{k,n}})_{n \in  [\![1,N]\!]}}$. Thanks to the majorizing condition \eqref{Eq_PropMaj}, the sequence $(x^k)_{k \in \mathbb{N}}$ defined in BD3MG algorithm satisfies, for every $k \in\mathbb{N}$, for every $c \in [\![1,C]\!]$,
\begin{center}
    \scalebox{0.90}[0.90]{$f(x^{k+1}) \leq Q_{(\mathcal{S}_c^k)}(x^{k+1}_{(\mathcal{S}_c^k)},x^{k,\delta_k}) \leq Q_{(\mathcal{S}_c^k)}(x^{k,\delta_k}_{(\mathcal{S}_c^k)},x^{k,\delta_k}) = f(x^{k,\delta_k}).$}
\end{center}



When there is no delay , i.e. $\delta_k \equiv 0$, we go back to a synchronous algorithm, benefiting from the classical block descent result inherent to MM schemes:
\begin{center}
     \scalebox{0.95}[0.95]{$f(x^{k+1}) \leq Q_{(\mathcal{S}_c^k)}(x^{k+1}_{(\mathcal{S}_c^k)},x^k) \leq Q_{(\mathcal{S}_c^k)}(x^{k}_{(\mathcal{S}_c^k)},x^k) = f(x^k)$.}
\end{center}

\subsection{Hypothesis and convergence}
\label{sec:HYPO}

We now state our convergence result for the sequence $(x^k)_{k\in \mathbb{N}}$ resulting from the BD3MG algorithm presented in Section \ref{alg:BMMD_master}. We first introduce the following assumptions.

\begin{assumption}
\label{Ass:A}
Function $f$ is differentiable, bounded from below and semi-algebraic\footnote{Real semi-algebraic functions represent a wide class of functions that satisfy the Kurdyka-\L{}ojasiewicz inequality, which is at the core of our convergence study. See \cite{Attouch2013,Bolte2014} for further details.}. Moreover, $f$ has an $\mathcal{L}$-Lipschitzian gradient on $\mathbb{R}^N$ with $\mathcal{L} > 0$, i.e.
    \vspace{0.1cm}
        \begin{equation}
            \forall (x,y) \in (\mathbb{R}^{N})^2, \; \|\nabla f(x)-\nabla f(y)\| \leq  \mathcal{L} \|x - y\|.
        \end{equation}
\end{assumption}

\begin{assumption}[Boundedness of delay]
\label{Ass:B}
Under the convention that $\forall \; l \in\mathbb{N}^*,\; \mathbb{S}_{-l} = \emptyset$, there exists $\tau \in \mathbb{N}$ such that
        \begin{equation}
             \forall k \in \mathbb{N}, \quad [\![1,N]\!]  \subset \bigcup\limits_{i=k-\tau}^{k} \mathbb{S}_i.
        \end{equation}
 \end{assumption}       
        
        \begin{assumption}[Curvature of quadratic majorant]
        \label{Ass:C}
Let us denote $(\Gamma^k_c)_{k\in\mathbb{N},c\in[\![1,C]\!] }$ the sequence of matrices defined as:
    \begin{equation*}
         \forall k\in \mathbb{N},\; \forall c \in [\![1,C]\!],\quad \Gamma^k_c = \mathcal{A}_{(\mathcal{S}^k_c)}(x^{k,\delta_k})-\frac{1}{2}\mathcal{A}_{(\mathcal{S}^k_c)}(x^{k}).
    \end{equation*}
        There exists $(\underline{\nu},\overline{\nu}) >0$ such that, for every $k\in \mathbb{N}$, for every $c \in [\![1,C]\!]$, \begin{equation}
        (\mathcal{L}\sqrt{\tau} + \underline{\nu}) \mathrm{Id}_{|\mathcal{S}^k_c| } \preceq \Gamma^k_c \preceq \overline{\nu} \mathrm{Id}_{|\mathcal{S}^k_c|}.
        \end{equation}        
        \end{assumption}

Under those assumptions, we state our convergence result for BD3MG, whose proof, relying on the analysis from \cite{davis2016asynchronous,Chouzenoux_2013}, is skipped by lack of space.
\begin{theorem}[Convergence of BD3MG] Let $f : \mathbb{R
}^N \mapsto \mathbb{R}$, and $x^0 \in \mathbb{R
}^N$. If Assumptions~\ref{Ass:A}-\ref{Ass:B}-\ref{Ass:C} are verified, then the sequence $(x^k)_{k\in\mathbb{N}}$ built by BD3MG algorithm converges globally to a stationary point $x^*$ of $f$.
\label{theo:BD3MG}
\end{theorem}
It is worthy to point out that, in contrast with most existing works in the literature of distributed optimization, no convexity assumption is required in our analysis. The MM framework allows to make use of the recent theory of non-smooth analysis \cite{Attouch2013}, as we also have shown in our previous works~\cite{Chouzenoux15jogo,Chouzenoux_2013}.

\input{Implementation.tex}

\small
\bibliographystyle{IEEEbib}
\bibliography{main}

\end{document}

%% file: Implementation.tex
\section{Application to 3D image deblurring}

\subsection{Model and objective function}
\begin{figure}[h]
    \centering
\begin{tabular}{@{}c@{}}
 \includegraphics[width=8cm]{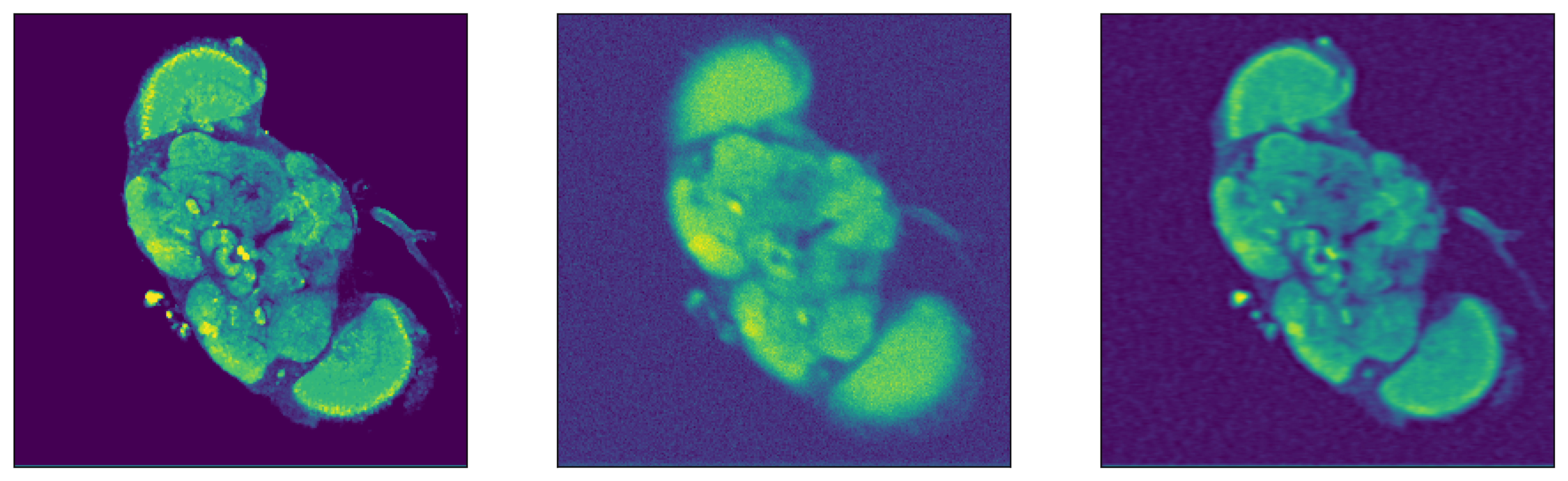}\\
    \includegraphics[width=8cm]{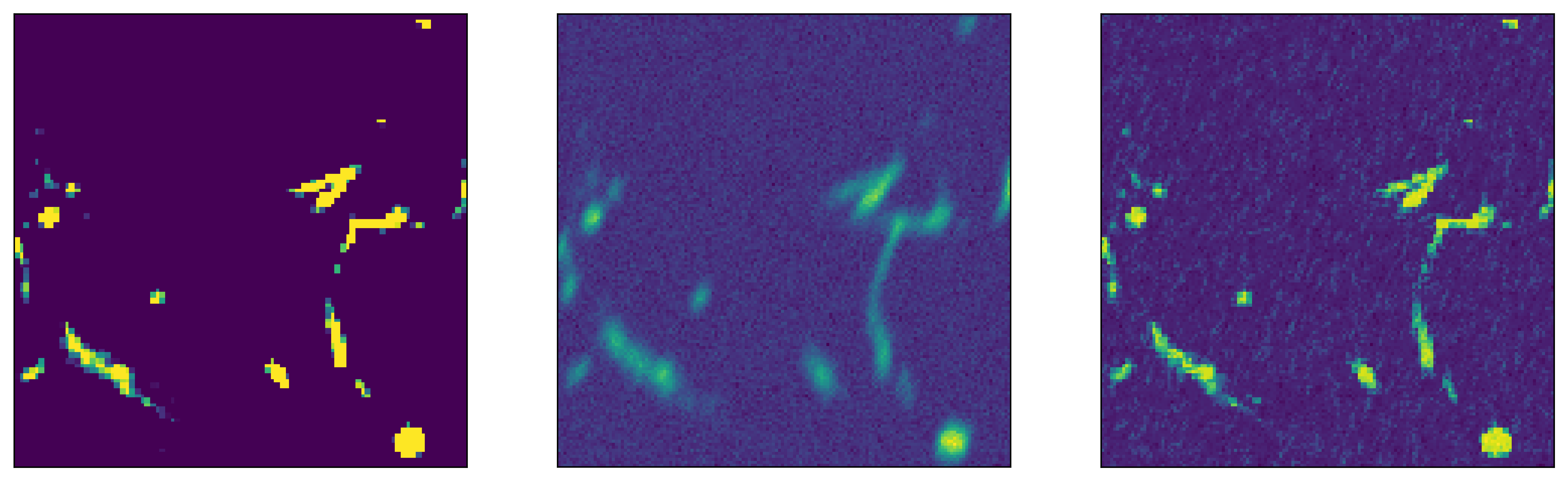}
\end{tabular}
    \caption{Comparison between original (\textbf{left}), degraded (\textbf{middle}) and restored (\textbf{right}) slices ($z=10$) of \texttt{FlyBrain} and \texttt{Aneurysm}.}
    \label{fig:visual}
\end{figure}

We consider the restoration of a 3D microscopic volume $\overline{x}$ of size $N = N_{\mathsf{X}} \times N_{\mathsf{Y}} \times N_{\mathsf{Z}}$, given a degraded observation $y$, altered by a depth-variant blur operator $H$ and an additive white noise $b$: 
\begin{equation}\label{eq:model}
y = H \overline{x} + b.
\end{equation}
The associated inverse problem can be solved efficiently by minimizing a least squares penalization penalized by a smooth 3D regularization function. Specifically, $f$ takes the form:
\begin{equation}\label{eq_crit}
\forall x \in \mathbb{R}^N, \quad f(x) = \sum_{s=1}^4 f_{s} (L_sx)
\end{equation} with $f_1 \circ L_1 = \frac{1}{2} \| H \cdot - y\|^2$, $f_2 \circ L_2 = \eta\, d^2_{[x_{\min},x_{\max}]^N}$,  $f_3 \circ L_3 = \lambda \sum_{n=1}^N \sqrt{ \left([V^{\mathsf{X}} \cdot]_n\right)^2 + \left( [V^{\mathsf{Y}} \cdot]_n\right)^2 + \delta^2}$, and $f_4 \circ L_4 = \kappa \| V^{\mathsf{Z}} \cdot \|^2$. 
Hereabove, $V^{\mathsf{X}} \in \eR^{N \times N}$, $V^{\mathsf{Y}} \in \eR^{N \times N}$, $V^{\mathsf{Z}} \in \eR^{N \times N}$ state for discrete gradient operators along the three directions of the volume, $x_{\min}$ (resp $x_{\max}$) $\in \eR$ are minimal (resp. maximal) bounds on the sought intensity values and $d_E$ is the Euclidian distance to set $E$. Function $f$ satisfies Ass.~\ref{Ass:A} as a combination of squared distances and discretized total-variation distances.  If not specified otherwise, the majorizing matrices for $f$ will be constructed following the strategy in \cite{Chouzenoux_2013}, ensuring the fulfillment of Ass.~\ref{Ass:C}. The convolution operator $H$ simulates a depth-varying 3D Gaussian blur. For each depth $z \in \{1,\ldots, N_{\mathsf{Z}}\}$, the blur kernel is characterized by different variance and rotation parameters whose values are taken as random realizations of uniform distributions. The depth-variant structure of the blur model motivates us to split the vector $x$ along the dimension $\mathsf{Z}$, assigning computing agents (i.e. workers) to one (or several) selected slice(s) of the 3D volume. For the sake of simplicity, a sequential ordering is used, for the selection rule of the processed blocks, hence Ass.~\ref{Ass:B} is fulfilled. Furthermore, the same strategy as in \cite{7532949} is employed to control the memory usage during communications master/worker. All algorithms are initialized with $x^0 \in \mathbb{R}^N$ whose entries are uniformly sampled in $[0,\max(y)]$. Parameters $\lambda, \eta, \delta, \kappa>0$ will be tuned manually, for both presented examples, so as to maximize the Signal-to-Noise Ratio (SNR) of the restored volume. Two microscopic images, namely \texttt{FlyBrain} and \texttt{Aneurysm} will be considered. Example of restoration results are presented in Fig.~\ref{fig:visual}.

\subsection{Ablation study and comparative analysis}

In order to measure the performance improvement allowed by BD3MG, we conducted an ``ablation study'', that consists in removing some/all acceleration features of our algorithm, namely asynchrony, subspace line-search and MM scaling. We then obtain, in addition to BD3MG, four parallel/distributed optimization algorithms that we detail hereafter. The \emph{asynchronous gradient descent} (Async-GD) corresponds to the algorithm from \cite{niu2011hogwild}, obtained by limiting the subspace to the sole gradient descent direction $\mathcal{D}_{(\mathcal{S})}(x) = - \nabla_{(\mathcal{S})}f(x)$ and by setting $\mathcal{A}_{(\mathcal{S})} = \mathcal{L} \, \mathrm{Id}_{|\mathcal{S}|}$, with $\mathcal{L}$ the Lipschitz constant of $\nabla f$. The \emph{asynchronous conjugate gradient} algorithm (Async-CG) identifies with our BD3MG method when using the basic MM metric $\mathcal{A}_{(\mathcal{S})} = \mathcal{L} \, \mathrm{Id}_{|\mathcal{S}|}$. The resulting algorithm can be viewed as a distributed version of the nonlinear conjugate gradient method, with closed form stepsize from \cite{Labat2008}. The \emph{asynchronous MM algorithm} (Async-MM) is obtained by removing the subspace acceleration strategy in BD3MG, so that $d'_{S} = - \mathcal{A}_{(\mathcal{S})}(x)^{-1} \nabla_{(\mathcal{S})}f(x)$ in the Worker loop. The latter inversion is performed with linear biconjugate gradient solver. Async-MM can be interpreted as a distributed implementation of a half-quadratic algorithm \cite{Allain2006}. Finally, \emph{block parallel 3MG algorithm} (BP3MG) is the method we originally proposed in  \cite{7532949}, where blocks updates are assumed to be computed at the same time without any communication delay. The names and characteristics of the resulting five tested methods are summarized in Tab.~\ref{tab:optisetup}.

\begin{table}[h]
  \centering
  \small
  \begin{tabular}{|c||ccc|} 
    \hline
    Name & Asynchrony & Memory & MM scaling \\
    \hline\hline
    Async-GD & \checkmark  & \xmark & \xmark  \\
    \hline
    Async-CG & \checkmark & \checkmark & \xmark   \\
    \hline
    Async-MM & \checkmark & \xmark & \checkmark \\
    \hline
    BP3MG & \xmark & \checkmark & \checkmark \\
    \hline
    \hline
    BD3MG & \checkmark & \checkmark & \checkmark\\
    \hline
  \end{tabular}
   \caption{Algorithmic features of the compared approaches.}
  \label{tab:optisetup}
\end{table} 
\vspace{-0.6cm}
\begin{figure}[h]
\includegraphics[width=8.5cm]{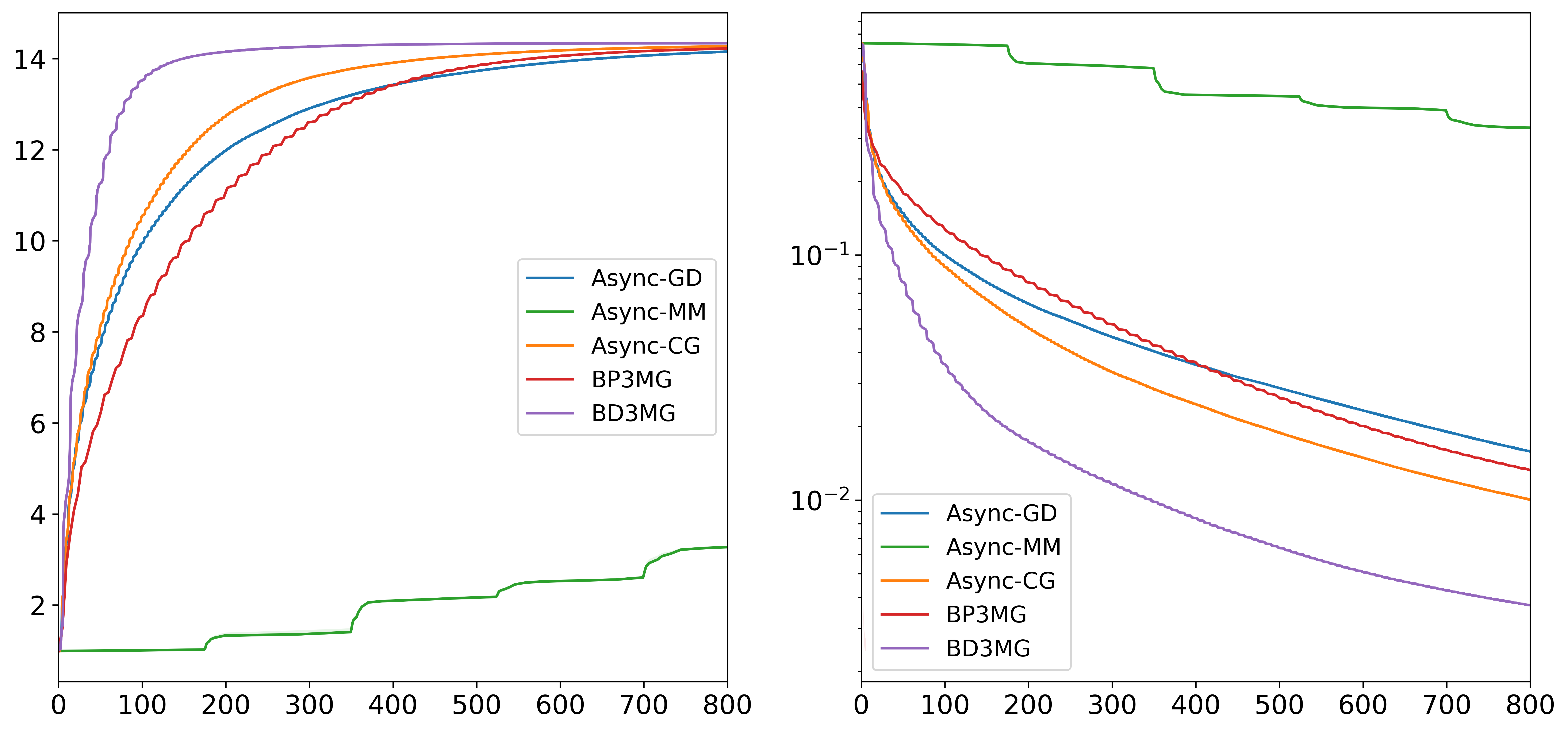}
\caption{Evolution of SNR in dB (left) and relative distance to solution $\|x^k-x^{*}\|/\|x^{*}\|$ (right) along time (in seconds) for \texttt{FlyBrain} restoration. Results are averaged over ten noise realizations.}
\label{fig:comp}
\end{figure}

%
%

We perform our comparative analysis, using the microscopic image  \texttt{FlyBrain} with size $N = 256 \times 256 \times 24$. It is degraded by depth-variant blur kernels of size $11\times 11 \times 21$, and a zero-mean white Gaussian noise with standard deviation $0.04$. The results are averaged over ten noise realizations, and the initial signal to noise ratio is around $11.6$ dB. Fig.~\ref{fig:comp} illustrates the evolution of the SNR of the restored image along time in seconds, for each method, for experiments performed using Python 3, running an Intel® Xeon(R) W-2135 CPU with 12 cores clocked at 3.70GHz. We also display the relative distance to the solution $x^*$, computed after a large number of iterations (typically, $10^4$), characterized by an average SNR of $14.33$ dB. One can see that the proposed BD3MG method clearly outperforms the others in terms of time to reach close-to-optimal solution. Async-MM led to the slowest convergence, as it requires linear system inversion at each iteration. The superiority of Async-CG over Async-GD illustrates the benefits for including the memory term within the subspace. Finally, BD3MG reaches convergence faster than its synchronous counterpart, BP3MG, probably thanks to the removal of the locking conditions in the Master loop. 

\subsection{Linear Speedup}

\begin{figure}[htb!]
\setlength{\belowcaptionskip}{-10pt}
\centering
\includegraphics[width=5cm]{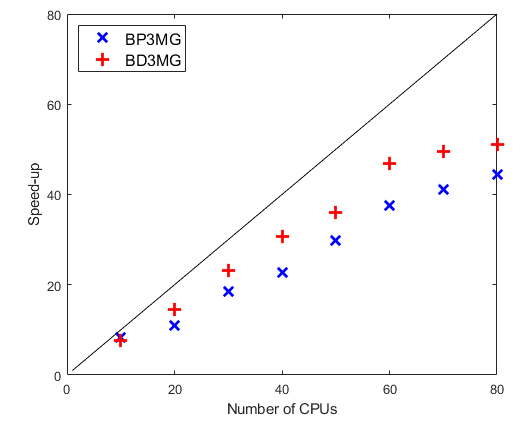}
\caption{Speed-up ratio for BD3MG (blue) and BP3MG (red), with respect to the number of active cores for the restoration of \texttt{Aneurysm}.}
\label{fig:speedup}
\end{figure}

In order to assess the scalability of the BD3MG algorithm, we further analysed the speed-up of the optimization process when a High Parallel Computing computer is being used. Namely, we ran BD3MG and BP3MG on an Intel Xeon CPU 6148 with up to 80 physical cores at 2.4 GHz (Skylake) and 1.5 Tio of RAM. Image \texttt{Aneurysm} with size $N = 155 \times 154 \times 79$ is degraded by blur kernels of size to $5 \times 5 \times 11$, and noise standard deviation of $0.04$, so that the initial SNR is $6.44$ dB, while the restored SNR is $11.92$ dB. Fig.~\ref{fig:speedup} presents the acceleration ratio between the required computation time for two cores (i.e. one Master and one Worker) versus the computation time when activating from 10 to 80 cores, for reaching the stopping criterion $\|x^{k+1}-x^k\| \leq 10^{-6} \|x^k\|$. This illustrates the great potential of scalability of the proposed algorithm. Asynchronicity of BD3MG allows to improve the speed-up, in comparison to the one exhibited by BP3MG \cite{7532949}. Finally, as the number of core increases, a mild saturation effect is observed (in
agreement with Amdahl’s law \cite{Amdahl}).

\section{Conclusion}\label{se:conclu}

In this paper, we have presented a new block distributed Majorize-Minimize Memory Gradient algorithm to tackle a wide class of large scale optimization problems. The main feature of our method lies in its distributed asynchronous formulation that allows for delays among workers, while theoretically maintaining the convergence guarantees and the practical performance of the powerful 3MG scheme. The new algorithm has been tested in the context of 3D image restoration under depth-variant blur. Experimental results underlined its scalability and efficiency\footnote{The code of BD3MG for depth-variant image deblurring has been made available at \texttt{https://github.com/mathieuchal/BD3MG}}. Future works will be focused on the extension to more general graph topologies.

%% file: main.bbl
\begin{thebibliography}{10}

\bibitem{Chouzenoux_2013}
E.~Chouzenoux, A.~Jezierska, J.-C. Pesquet, and H.~Talbot,
\newblock ``A majorize-minimize subspace approach for $\ell_2-\ell_0$ image
  regularization,''
\newblock {\em SIAM Journal on Imaging Sciences (SIIMS)}, vol. 6, no. 1, pp.
  563–591, January 2013.

\bibitem{Tseng2001}
P.~Tseng,
\newblock ``Convergence of a block coordinate descent method for
  nondifferentiable minimization,''
\newblock {\em Journal of Optimization, Theory and Applications}, vol. 109, no.
  3, pp. 475--494, 2001.

\bibitem{4303144}
M.~W. {Jacobson} and J.~A. {Fessler},
\newblock ``An expanded theoretical treatment of iteration-dependent
  majorize-minimize algorithms,''
\newblock {\em IEEE Transactions on Image Processing}, vol. 16, no. 10, pp.
  2411--2422, Oct 2007.

\bibitem{Hong2016}
M.~{Hong}, M.~{Razaviyayn}, Z.~{Luo}, and J.~{Pang},
\newblock ``A unified algorithmic framework for block-structured optimization
  involving big data: With applications in machine learning and signal
  processing,''
\newblock {\em IEEE Signal Processing Magazine}, vol. 33, no. 1, pp. 57--77,
  Jan 2016.

\bibitem{Bolte2014}
J.~Bolte, S.~Sabach, and M.~Teboulle,
\newblock ``Proximal alternating linearized minimization for nonconvex and
  nonsmooth problems,''
\newblock {\em Mathematical Programming}, vol. 146, no. 1, pp. 459--494, Aug
  2014.

\bibitem{LeeNMF}
D.~D. Lee and H.~S. Seung,
\newblock ``Algorithms for non-negative matrix factorization,''
\newblock in {\em Proceedings of the 13th International Conference on Neural
  Information Processing Systems (NIPS 2000)}, Denver, Colorado, 2000, p.
  535–541.

\bibitem{Fessler97}
J.~A. Fessler,
\newblock ``{Grouped coordinate descent algorithms for robust edge-preserving
  image restoration},''
\newblock in {\em Image Reconstruction and Restoration II}, Timothy~J. Schulz,
  Ed. International Society for Optics and Photonics, 1997, vol. 3170, pp. 184
  -- 194, SPIE.

\bibitem{Repetti2015}
L.~Duval E.~Chouzenoux A.~Repetti, M. Q.~Pham and J.-C. Pesquet,
\newblock ``Euclid in a taxicab: Sparse blind deconvolution with smoothed
  $\ell_1$/$\ell_2$ regularization,''
\newblock {\em IEEE Signal Processing Letters}, vol. 22, no. 5, pp. 539--543,
  May 2015.

\bibitem{Chouzenoux15jogo}
E.~Chouzenoux, J.-C. Pesquet, and A.~Repetti,
\newblock ``A block coordinate variable metric forward-backward algorithm,''
\newblock {\em Journal of Global Optimization}, vol. 66, no. 3, pp. 457--485,
  2015.

\bibitem{Fessler2002}
S.~Sotthivirat and J.~A. Fessler,
\newblock ``Image recovery using partitioned-separable paraboloidal surrogate
  coordinate ascent algorithms,''
\newblock {\em {IEEE} {T}ransactions on {S}ignal {P}rocessing}, vol. 11, no. 3,
  pp. 306--317, 2002.

\bibitem{7532949}
S.~{Cadoni}, E.~{Chouzenoux}, J.~{Pesquet}, and C.~{Chaux},
\newblock ``A block parallel majorize-minimize memory gradient algorithm,''
\newblock in {\em Proceedings of the IEEE International Conference on Image
  Processing (ICIP 2016)}, Phoenix, Arizona, 25-28 Sep. 2016, pp. 3194--3198.

\bibitem{Moon2019}
G.~E. Moon, A.~Sukumaran{-}Rajam, S.~Parthasarathy, and P.~Sadayappan,
\newblock ``{PL-NMF:} parallel locality-optimized non-negative matrix
  factorization,'' 2019,
\newblock http://arxiv.org/abs/1904.07935.

\bibitem{niu2011hogwild}
F.~Niu, B.~Recht, C.~Re, and S.~J. Wright,
\newblock ``Hogwild: A lock-free approach to parallelizing stochastic gradient
  descent,''
\newblock in {\em Proceedings of the 25th Conference on Advances in Neural
  Information Processing Systems (NIPS 2011)}, pp. 693--701. Granada, Spain,
  12-17 Dec. 2011.

\bibitem{grishchenko2018asynchronous}
D.~Grishchenko, F.~Iutzeler, J.~Malick, and M.R. Amini,
\newblock ``Asynchronous distributed learning with sparse communications and
  identification,'' 2018,
\newblock https://arxiv.org/abs/1812.03871.

\bibitem{Zhang:2014:ADA:3044805.3045082}
R.~Zhang and J.~T. Kwok,
\newblock ``Asynchronous distributed {ADMM} for consensus optimization,''
\newblock in {\em Proceedings of the 31st International Conference on Machine
  Learning (ICML 2014)}, 21-26 June 2014, pp. 1701--1709.

\bibitem{pesquet2014class}
J.-C. Pesquet and A.~Repetti,
\newblock ``A class of randomized primal-dual algorithms for distributed
  optimization,''
\newblock {\em Journal on Nonlinear Convex Analysis}, vol. 16, no. 12, pp.
  2353--2490, Dec. 2015.

\bibitem{Hannah_2017}
R.~Hannah and W.~Yin,
\newblock ``On unbounded delays in asynchronous parallel fixed-point
  algorithms,''
\newblock {\em Journal of Scientific Computing}, vol. 76, no. 1, pp. 299–326,
  Dec 2017.

\bibitem{Abboud2019}
F.~Abboud, E.~Chouzenoux, J.-C. Pesquet, and H.~Talbot,
\newblock ``Distributed algorithms for proximity operator computation with
  applications to video processing,'' 2019,
\newblock https://hal.archives-ouvertes.fr/hal-01942710.

\bibitem{Onose16}
A.~Onose, R.~E. Carrillo, A.~Repetti, J.~D. McEwen, J.-T. Thiran, J.-C.
  Pesquet, and Y.~Wiaux,
\newblock ``Scalable splitting algorithms for big-data interferometric imaging
  in the {SKA} era,''
\newblock {\em Monthly Notices of the Royal Astronomical Society}, vol. 462,
  no. 4, pp. 4314--4335, 2016.

\bibitem{Xu2019}
J.~Xu, Y.~Sun, Y.~Tian, and G.~Scutari,
\newblock ``A unified contraction analysis of a class of distributed algorithms
  for composite optimization,'' 2019,
\newblock https://arxiv.org/abs/1910.09817.

\bibitem{DavisNIPS}
D.~Davis, M.~Udell, and B.~Edmunds,
\newblock ``The sound of {APALM} clapping: Faster nonsmooth nonconvex
  optimization with stochastic asynchronous palm,''
\newblock in {\em Proceedings of the 30th International Conference on Neural
  Information Processing Systems (NIPS 2016)}, 5-10 Dec. 2016, p. 226–234.

\bibitem{davis2016asynchronous}
D.~Davis,
\newblock ``The asynchronous palm algorithm for nonsmooth nonconvex problems,''
  2016,
\newblock https://arxiv.org/abs/1604.00526.

\bibitem{Thouvenin}
P.~{Thouvenin}, N.~{Dobigeon}, and J.~{Tourneret},
\newblock ``Partially asynchronous distributed unmixing of hyperspectral
  images,''
\newblock {\em IEEE Transactions on Geoscience and Remote Sensing}, vol. 57,
  no. 4, Apr. 2019.

\bibitem{Sun2016}
Y.~Sun, P.~Babu, and D.~P. Palomar,
\newblock ``Majorization-minimization algorithms in signal processing,
  communications, and machine learning,''
\newblock {\em IEEE Transactions on Signal Processing}, vol. 65, no. 3, pp.
  794--816, 2016.

\bibitem{Attouch2013}
H.~Attouch, J.~Bolte, and B.F. Svaiter,
\newblock ``Convergence of descent methods for semi-algebraic and tame
  problems: proximal algorithms, forward--backward splitting, and regularized
  gauss--seidel methods,''
\newblock {\em Mathematical Programming}, vol. 137, no. 1, pp. 91--129, Feb
  2013.

\bibitem{Labat2008}
C.~Labat and J.~Idier,
\newblock ``Convergence of conjugate gradient methods with a closed-form
  stepsize formula,''
\newblock {\em Journal of Optimization Theory and Applications}, vol. 136, no.
  1, pp. 43--60, Jan 2008.

\bibitem{Allain2006}
M.~{Allain}, J.~{Idier}, and Y.~{Goussard},
\newblock ``On global and local convergence of half-quadratic algorithms,''
\newblock {\em IEEE Transactions on Image Processing}, vol. 15, no. 5, pp.
  1130--1142, May 2006.

\bibitem{Amdahl}
G.~M. Amdahl,
\newblock ``Validity of the single processor approach to achieving large scale
  computing capabilities,''
\newblock in {\em Proceedings of the American Federation of Information
  processing Societies conference (AFIPS 1967)}, Atlantic City, 18-20 Apr.
  1967, pp. 483--485.

\end{thebibliography}
